\documentclass[11pt]{amsart}
\usepackage{latexsym, amssymb, amsmath, amscd,mathrsfs}
 \usepackage{eucal}

\title{A Local to Global Question for Linear Functionals}

   \author{George F. Seelinger and Wenhua Zhao$^\dag$}  
   \thanks{$^\dag$ Wenhua Zhao passed away in December of 2023.  We were in the process of writing this paper when he passed.}    
    \date{\today}
 \address{G. F. Seelinger, Department of Mathematics, Illinois State University, Normal, IL 61761. {\it Email}: gfseeli@ilstu.edu}
\address{W. Zhao, Department of Mathematics, Illinois State University, Normal, IL 61761.}
 
\date{\today}

\newcommand{\A}{{\mathbb{A}}}

\newcommand{\uu}{{\vec{u}}}
\newcommand{\vv}{{\vec{v}}}
\newcommand{\ww}{{\vec{w}}}

\newcommand{\yy}{{\vec{y}}}
\newcommand{\zz}{{\vec{z}}}

\newcommand{\vsp}{\vspace{1em}}

\newcommand{\spn}{{\mathrm{Span}}}

\newcommand{\prf}{{\em Proof.} \ \ }

\newcommand{\spec}{{\mathrm{Spec}}}

\newcommand{\chr}{{\mathrm{char}}}

\newcommand{\ds}{\displaystyle}

\newcommand{\ve}{\vec{e}}

\newcommand{\vq}{\vec{q}}
\newcommand{\veca}{\vec{a}}
\newcommand{\vecb}{\vec{b}}
\newcommand{\vecp}{\vec{p}}
\newcommand{\vecpp}{\vec{p}^{\,'}}
\newcommand{\vecap}{\vec{a}^{\,'}}
\newcommand{\vecbp}{\vec{b}^{\,'}}

\newcommand{\veccp}{\vec{c}^{\,'}}

\newtheorem{thm}{Theorem}

\newtheorem{lemma}[thm]{Lemma}
\newtheorem{cor}[thm]{Corollary}

\begin{document}

\maketitle

\begin{abstract}
Let $F$ be an algebraically closed field and let $n\geq 3$.  Consider 
$V=F^n$ with standard basis
$\{\ve_1,\ldots,\ve_n\}$ and its dual space $V^*=
{\mathrm{Hom}}_{F-{\mathrm{lin}}}(V,F)$ with dual basis $\{y_1,\ldots,y_n\}\subseteq
V^*$ and let $\yy = \sum_i y_i\otimes e_i\in V^*\otimes V$.  Let $d<n$ and consider the vectors $\vq_1,\ldots,\vq_d\in V^*\otimes V$.   In this note we consider the question of 
whether $\yy(\vv)=\vv \in \spn_F(\vq_1(\vv),\ldots,\vq_d(\vv))$ for all $\vv\in V$ implies
that $\yy\in \spn_F(\vq_1,\ldots,\vq_d)$.  We show this is true for $d=1$ or $d=2$, but 
that additional properties are needed for $d\geq 3$.  We then interpret this result in terms
of subspaces of $M_n(F)$ that do not contain any rank 1 idempotents.
\end{abstract}

\section{Introduction}

Let $F$ be an algebraically closed field and let $n>2$ be a positive integer.  Assume 
that ${\mathrm{char}}(F)>n$ or $\chr(F)=0$ and let $V=F^n$.  Let $V^*={\mathrm{Hom}}_{F-{\mathrm{lin}}}(V,F)$
be the space of linear functionals with basis $\{y_1,\ldots,y_n\}$ where $y_i(\vv)=v_i$ is just projection onto the $i$th component.
Let \[ \yy = \left[ \begin{array}{c} y_1 \\ y_2 \\ \vdots \\ y_n
\end{array}\right]\in (V^*)^n\] and $\{\vq_1,\ldots,\vq_d\}$ be an $F$-linearly independent set in $(V^*)^n$ with $d<n$ such that
for all $\vv\in V$, $\yy$ evaluated at $\vv$ is in the span of the vectors $\vq_i$ evaluated at $\vv$.  In this note we ask
under what conditions can we conclude that $\yy\in \spn_F(\vq_1,\ldots,\vq_d)$?

We can make the above question more precise as follows.  Let $S=F[y_1,\ldots,y_n]$ be the polynomial ring in $n$ commuting variables
and let $S_1 = \spn_F(y_1,\ldots,y_n)\subseteq S$ be the subspace of elements that are homogeneous of degree $1$.  Hence we can identify 
$V^*$ with $S_1$.  For any maximal ideal ${\frak m}$, let $\phi_{{\frak m}}:S^n\rightarrow (S/{\frak m})^n$ be the natural projection
defined by letting $\phi_{{\frak m}}(\ww)$ be the vector whose $i$th component is $w_i+{\frak m}\in S/{\frak m}$ for all $i$ for any $\ww\in S^n$.
Let $\yy,\vq_1,\ldots,\vq_d$ be as above.  Then we can ask if $\phi_{{\frak m}}(\yy)\in \spn_F(\phi_{{\frak m}}(\vq_1),\ldots,\phi_{{\frak m}}(\vq_d))$
for all maximal ideals ${\frak m}\subseteq S$, is $\yy\in \spn_F(\vq_1,\ldots,\vq_d)$?  In this note we show that without additional conditions, this
is only true if $d=1$ or $d=2$.

While we believe the above question is of general interest, 
our particular interest in this question arose from our study of Mathieu-Zhao subspaces of matrix algebras.  In \cite{Zhao1}, the second author first introduced the concept of a Mathieu subspace.  
Subsequently, in the literature, these Mathieu subspaces started being called Mathieu-Zhao subspaces or Mathieu-Zhao spaces (e.g., see \cite{EKC, DEZ}).  The introduction of Mathieu-Zhao subspaces was directly motivated by the Mathieu Conjecture \cite{Ma} and the Image Conjecture \cite{IC}, each of which implies the Jacobian conjecture \cite{K, BCW, E}. For example, the Jacobian conjecture will follow if some explicitly given subspaces of multivariate polynomial algebras over $\mathbb C$ can be shown to be 
Mathieu-Zhao subspaces of these polynomial algebras. For details, see \cite{IC, DEZ, EKC}. 

Note that ideals of rings are Mathieu-Zhao subspaces, but not conversely. Therefore  the concept of 
a Mathieu-Zhao subspace can be viewed as a natural generalization of the concept of an ideal. However, in contrast to ideals, Mathieu-Zhao subspaces are currently far from being well-understood. This is even the case for the most of finite rings or finite dimensional algebras over a field. For example, the classification of all Mathieu-Zhao subspaces 
of the matrix algebras $M_n(F)$ when $n\ge 3$ is still wide open.

The study of Mathieu-Zhao subspaces of matrix algebras $M_n(F)$ was initiated by the second author in 
\cite{Zhao2} in which the following two results were proven.  

\begin{thm}{\cite[Thm.\ 4.2]{Zhao2}}\label{ZIT}  Let $V\subseteq M_n(F)$ be a proper subspace of $M_n(F)$.  Then $V$ is a Mathieu-Zhao subspace of $M_n(F)$ if and only if $V$ does not contain any nonzero idempotents of $M_n(F)$.
\end{thm}

\begin{thm}\cite[Thm $5.1$]{Zhao2}\label{MS5.1}
Let $H$ be the subspace of $M_n(F)$ consisting of all trace-zero matrices. Then the following statements holds:
\begin{enumerate}
  \item[$i)$] if char.\,$F=0$ or char.\,$F>n$, then $H$ is the only Mathieu-Zhao subspace of $M_n(F)$ of codimension $1$;  
  \item[$ii)$] if char.\,$F\le n$, then $M_n(F)$ has no Mathieu-Zhao subspaces of  codimension $1$. 
\end{enumerate}
\end{thm}

Actually, the theorem above holds also for one-sided codimension one 
Mathieu-Zhao subspaces of $M_n(F)$. 

We note that given the characterization of Mathieu-Zhao subspaces of $M_n(F)$ given in Theorem \ref{ZIT}, for the purposes of this note we
use this as our definition of a Mathieu-Zhao subspace of $M_n(F)$ instead of stating the original definition here.

It is easy to see from the two theorems above that every subspace of a 
Mathieu-Zhao subspace of $M_n(F)$ is also a Mathieu-Zhao subspace of $M_n(F)$. Therefore, to classify all Mathieu-Zhao subspaces of $M_n(F)$ it suffices to classify all maximal 
Mathieu-Zhao subspaces of $M_n(F)$. In 2016, M. deBondt proved that the only Mathieu-Zhao
subspaces of $M_n(F)$ of codimension $1\leq d<n$ are subspaces of $H$ (see \cite{deBondt}).  

In our study of Mathieu-Zhao subspaces of $M_n(F)$, we asked whether we could characterize subspaces of $M_n(F)$ of
codimension $1\leq d <n$ that do not contain any rank 1 idempotents. We will say a subspace $W\subseteq M_n(F)$
is {\em r1-free subspace} if $W$ does not contain any rank 1 idempotents.  As part of our effort to characterize r1-free subspaces
of codimension $1\leq d <n$, we ask if there are such subspaces that are not subspaces of the trace zero matrices.  In this
paper, we give an anaology to deBondt's Theorem for r1-free subspaces and give an example of a codimension 3 r1-free subspace that is not contained in the trace zero matrices 
for $n\geq 4$.

\section{The $\yy$-local membership property}

 For any $F$-subspace $V\subseteq (S_1)^n$, Let $S.V$ denote the $S$-submodule of $S^n$ generated by $V$.  So if $\{\vq_1,\ldots,\vq_d\}$ is an $F$-basis
of $V$, then 
\[ S.V = \left\{ s_1\vq_1 + \cdots + s_d\vq_d : s_1,\ldots,s_d\in S\right\}. \]
We further define 
\[ V_L =  \left\{ c_1\vq_1 + \cdots + c_d\vq_d : c_1,\ldots,c_d\in L\right\} \subseteq L^n. \]
We note that since $V_L$ is generated by $V$ as an $L$-vector space, we always have
$\dim_L(V_L)\leq \dim_F(V)$.  Furthermore, if $\dim_L(V_L)=\dim_F(V)$, then $S.V$ is a free $S$-module of rank $d=\dim_F(V)$ as any $L$-basis $\{\vq_1,\ldots,\vq_d\}\subseteq V$ is also a basis of $S.V$ as a free $S$-module.  So the subspaces $V\subseteq (S_1)^n$ such that $\dim_F(V)=
\dim_L(V_L)$ are minimal among the set of subspaces $V'\subseteq (S_1)^n$ such that
$(V')_L=V_L$.  

For any maximal ideal ${\frak m}\subseteq S$ and any
$z\in S$, we let $[z]_{{\frak m}}$ denote the image of $z$ in $\kappa({\frak m})=S/{\frak m}$.  Furthermore, 
we define 
$\phi_{{\frak m}}: S^n \rightarrow (S^n)\otimes_S \kappa({\frak m}) \cong \kappa({\frak m})^n$
be the induced projection on $S^n$.  (i.e., If $\vv\in S^n$, then $\phi_{{\frak m}}(\vv)$ just reduces each component of $\vv$ modulo ${\frak m}$.)
We note that for each $\vq\in (S_1)^n$, there exists a $b\in M_n(F)$ such that $\vq = b\yy$.  In this case, $\phi_{{\frak m}}(\vq) = b\,\phi_{{\frak m}}(\yy)$.

We are interested in the implications of the following local property:

\vsp

\noindent{{\bf $\yy$-Local Membership Property:}} 
Let $V\subseteq (S_1)^n$ be an $F$-subspace.
We say $V$ has the {\em $\yy$-local membership property} if 
for evey maximal ideal ${\frak m}\subseteq S$, we have $\phi_{{\frak m}}(\yy) \in
\phi_{{\frak m}}(S.V)$.

\vsp

In particular, the main question that motivated this note can be phrased as follows:

\vsp

\noindent{{\bf Main Question:}} Let $V\subseteq (S_1)^n$ be an $F$-subspace.
If $V$ has the $\yy$-local membership property, what other conditions are sufficient to 
imply $\yy\in V$?

\vsp

For $n\geq 4$ and $d\geq 3$, we give the following example of a subspace $V\subseteq S_1$
such that $d=\dim_F(V)=\dim_L(V_L)$ with the $\yy$-local membership property, but
$\yy\not\in V$. Hence extra conditions are needed to imply $\yy\in V$.  

\vsp

\noindent{{\bf Example:}}  Let $\{\ve_1,\ldots,\ve_n\}$ be the standard basis of 
$S^n$.  Let $\vq_1 = \yy -y_1\ve_{n-1}, \vq_2 = y_1\ve_{n-1}-y_2\ve_n, \vq_3 =
y_1\ve_n$ and let $\{\vq_4,\ldots,\vq_d\}$ be an $L$-linearly independent
subset of $\spn_S(\ve_1,\ldots,\ve_{n-3})$.  Hence $V=\spn_F(\vq_1,\ldots,\vq_d)$ is
such that $\dim_{L}(V_L)=d$ and $\yy\not\in\spn_L(\vq_4,\ldots,\vq_d)$.  

Now let ${\frak m}
\subseteq S$ be a maximal ideal.  If $y_1\in {\frak m}$, then $\phi_{{\frak m}}(\yy) =
\phi_{{\frak m}}(\vq_1)\in \phi_{{\frak m}}(S.V)$.  If $y_1\not\in
{\frak m}$, then $[y_1]_{{\frak m}}$ is invertible in $S/{\frak m}$.  So we can write
\[ \phi_{{\frak m}}(\yy) = \phi_{{\frak m}}(\vq_1) + \phi_{{\frak m}}(\vq_2) + 
[y_2]_{{\frak m}}[y_1]_{{\frak m}}^{-1}\phi_{{\frak m}}(\vq_3)\in 
\phi_{{\frak m}}(S.V). \] 
Hence $V$ has the $\yy$-local membership property.

Finally, if we set up the equation 
\[ \yy = \alpha_1\vq_1+\ldots + \alpha_d\vq_d \ \mbox{ where } \  \alpha_1,\ldots,
\alpha_d\in F, \]
we get 
\[ y_{n-1} = \alpha_1(y_{n-1}-y_1) + \alpha_2y_1 \ \mbox{ and } \ 
y_n = \alpha_1y_n-\alpha_2y_2 + \alpha_3y_1. \]
This system has no solution in $F$ as the second equation implies $\alpha_1=1,\alpha_2=
\alpha_3=0$, which is not a solution to our first equation.  Therefore, 
$\yy\not\in V$. 

\vsp

In Theorem \ref{d=1or2}, we show that when $\dim_L(V_L)=1$ or $\dim_{L}(V_L)=2$ with $d<n$ we have
$\yy\in V$ if and only if $V$ has the $\yy$-local membership property.  
In addition, we explore other sufficient conditions on $V$ for which the $\yy$-local 
membership property implies $\yy\in V$.  First, we make the following observation.

\begin{lemma}\label{Lspan}
Let $V\subseteq (S_1)^n$ be a $d$-dimensional subspace such that $d=\dim_L(V_L)\leq n$.
If $V$ has the $\yy$-local membership property then $\yy\in V_L$.
\end{lemma}

\prf  First note that if $d=n$, then $V_L=L^n$, so 
our conclusion holds.  So assume $d<n$ and let 
$\{\vq_1,\ldots,\vq_d\}\subseteq	 V$ be an $L$-basis of $V_L$.  
If $\yy\not\in V_L$, then there must be a nonzero $(d+1)\times (d+1)$ minor $\Delta\in S$ of the matrix 
$C_{\yy} = \left[\begin{array}{ccccc} | & | & & | & | \\
\vq_1 & \vq_2 & \cdots & \vq_d & \yy \\
| & | & & | & | \end{array}\right]$.  Let ${\frak m}\subseteq S$ be a maximal ideal  such that 
$\Delta\not\in{\frak m}$.  Then $[\Delta]_{{\frak m}}\neq 0\in S/{\frak m}$, so it follows that
$\phi_{{\frak m}}(\yy)\not\in \phi_{{\frak m}}(S.V)$.  Hence our lemma follows.  \qed

\vsp

Note that the converse of Lemma \ref{Lspan} does not hold in general.  Indeed, if $n\geq 3$, we let $\vq_1 = \yy - y_2\ve_2$ and $\vq_2=y_1\ve_1$, then
we get $\yy\in\spn_{L}(\vq_1,\vq_2)$, but for any maximal ideal ${\frak m}$ such that $y_1\in {\frak m}$ but $y_2\not\in{\frak m}$, we have
$\phi_{{\frak m}}(\yy)\not\in\spn_F(\phi_{{\frak m}}(\vq_1),\phi_{{\frak m}}(\vq_2))$.  

Before ending this section, we give characterization of the $\yy$-local membership property in
terms of a morphism of affine varieties.

Let $V\subseteq (S_1)^n$ be a subspace such that $S.V$ is a free $S$-submodule of 
rank $d=\dim_F(V)$.  Let $\{\vq_1,\ldots,\vq_d\}$ be a basis of $V$ that is also a basis
of $S.V$.  Then we can consider a closed subvariety $Z_V$ of 
$\A_F^{n+d} = \spec_F(S[c_1,\ldots,c_d])$ defined by the ideal 
\[ P_V = \langle c_1q_{1,j} + \cdots + c_dq_{d,j} - y_i \ : \ 1\leq j\leq n\rangle , \]
where $c_1,\ldots,c_d$ are commuting indeterminants. Hence 
\[ Z_V = \spec_F\left(S[c_1,\ldots,c_d]/P_V\right). \]
Let $\theta: Z_V\rightarrow \A^n_F$ be the morphism defined by the comorphism
$\theta^{\sharp}:S\rightarrow S[c_1,\ldots,c_d]/P_V$ where 
$\theta^{\sharp}$ is the unique $F$-algebra homomorphims such that
$\theta^{\sharp}(y_j)=y_j+P_V$ for all $1\leq j\leq n$.  Then we get the following lemma.

\begin{lemma}\label{Surjective}
Let $V\subseteq (S_1)^n$ be an $F$-subspace such that 
$S.V$ is a free $S$-module of rank $d=\dim(V)$.  Then $V$ has the 
$\yy$-local membership property if and only if there exists a basis
$\{\vq_1,\ldots,\vq_d\}$ of $V$ such that $\theta:Z_V\rightarrow \A^n_F$ is 
surjective.  
\end{lemma}

\prf  Assume $V$ has the $\yy$-local membership property and choose 
a basis $\{\vq_1,\ldots,\vq_d\}\subseteq V$.  Let ${\frak m}\in S$ be a maximal 
ideal and let $\theta:Z_V\rightarrow \A^n_F$ be the affine morphism defined above corresponding to 
our basis.  Then by the $\yy$-local membership property, we have
$\phi_{{\frak m}}(\yy)\in \spn_F(\phi_{{\frak m}}(\vq_1),\ldots,\phi_{{\frak m}}(\vq_d))$, 
so there exist $\alpha_1,\ldots,\alpha_d\in F$ such that
\[ \alpha_1\phi_{{\frak m}}(\vq_1)+\cdots+\alpha_d\phi_{{\frak m}}(\vq_d) - 
\phi_{{\frak m}}(\yy) = \vec{0}. \]
Then
\[ {\frak m}' =  \langle c_1-\alpha_1 ,\ldots, c_d-\alpha_d\rangle + S[c_1,\ldots,c_d]{\frak m}\]
is a maximal ideal of $S[c_1,\ldots,c_d]$ containing $P_V$, hence defines a maximal ideal
of $S[c_1,\ldots,c_d]/P_V$ such that $(\theta^{\sharp})^{-1}({\frak m}')={\frak m}$.  Therefore, $\theta$ is surjective.

Conversely, assume $\theta$ defined by the basis $\{\vq_1,\ldots,\vq_d\}$ of $V$ 
is surjective and let ${\frak m}$ be a maximal ideal of $S$.  Since $\theta$ is surjective,
there exists a maximal ideal ${\frak m}'\subseteq S[c_1,\ldots,c_d]/P_V$ such that
$(\theta^{\sharp})^{-1}({\frak m}')={\frak m}$.  Let ${\frak m}''\subseteq S[c_1,\ldots,
c_d]$ be a maximal ideal containing $P_V$ such that ${\frak m}''/P_V={\frak m}'$.  Since
$F$ is algebraically closed, we get ${\frak m}'=\langle c_1-\alpha_1,\ldots,
c_d-\alpha_d\rangle+{\frak m}S[c_1,\ldots,c_d]$ for some $\alpha_1,\ldots,\alpha_d\in F$.  Then it follows that
\[ \phi_{{\frak m}}(\yy)=\alpha_1\phi_{{\frak m}}(\vq_1)+\cdots 
+ \alpha_d\phi_{{\frak m}}(\vq_d), \]
so $\phi_{{\frak m}}(\yy)\in \phi_{{\frak m}}(V)$.  Since ${\frak m}$ was arbitrary, $V$ has
the $\yy$-local membership property.  \qed

\section{Some Local to Global Results}

Let $V\subseteq (S_1)^n$ be such that $S.V$ is a free
$S$-module of rank $d=\dim_F(V)$ and has the $\yy$-local membership membership 
property.  In this section, we explore when $\yy\in V$.  In particular, if $n> 2\geq d$,
then $\yy\in V$.  

Given Lemma \ref{Lspan}, we first prove some results when $V$ satisfies the weaker 
condition of $\yy\in V_L$.  In proving these results, we choose an $F$-basis
$\{ \vq_1,\ldots,\vq_d\}$ of $V$ that is also an basis of $V_L$.  We define 
the following matrices relative to our basis:
\[ Q = \left[ \begin{array}{cccc} | & | & & | \\
\vq_1 & \vq_2 & \cdots & \vq_d \\
| & | & & | \end{array}\right] \ \mbox{ and for any $\ww\in(S_1)^n$, we define }
C_{\ww} = \left[ \ Q \ | \ \ww\, \right]. \]
Furthermore, if $I=\{i_1<i_2<\cdots<i_k\}\subseteq \{1,\ldots,n\}$ for some $k\leq n$, 
we will let 
\[ Q_I = \left[\begin{array}{cccc}
q_{1,i_1} & q_{2,i_1} & \cdots & q_{d,i_1} \\
q_{1,i_2} & q_{2,i_2} & \cdots & q_{d,i_2} \\
\vdots & & \ddots & \vdots \\
q_{1,i_k} & q_{2,i_k} & \cdots & q_{d,i_k} \end{array}\right], \ \ 
(\ww)_I = \left[ \begin{array}{c} w_{i_1} \\ w_{i_2} \\ \vdots \\ w_{i_k} \end{array}
\right], \ \mbox{ and } \ (C_{\ww})_I = \left[\ Q_I \ | \ (\ww)_I \ \right]. \]
With this notation, we get the following theorem.

\begin{thm}\label{yLthm}
Let $V\subseteq (S_1)^n$ be such that $S.V$ is a free module of rank $d=\dim_F(V)$
and let $\{\vq_1,\ldots,\vq_d\}$ be an $F$-basis of $V$.  Then for any $\ww\in V_L\cap
(S_1)^n$, there exist homogeneous $h_1,\ldots,h_d,k_1,\ldots,k_d\in S$ such that
\[ \ww = \frac{h_1}{k_1}\vq_1 + \cdots + \frac{h_d}{k_d}\vq_d \]
where for all $j$, $\gcd(h_j,k_j)=1$ and if $h_j\neq 0$ then $\deg(h_j)=\deg(k_j)\leq d$.
Furthermore, if $m={\mathrm{lcm}}(\{k_j : 1\leq k\leq d\})$ then
$\det(Q_I)\in \langle m\rangle \subseteq S$ for all $I\subseteq \{1,\ldots,n\}$ of order $d$.
\end{thm}

\prf  Let $\{\vq_1,\ldots,\vq_d\}$ be an $F$-basis of $V$.  Then 
$\{\vq_1,\ldots, \vq_d\}$ is a basis of $V_L$ since $V_L$
is generated by $\{\vq_1,\ldots,\vq_d\}$ and any dependence relation over $L$
would imply a dependence relation over $S$, contradicting that $S.V$ is a free $S$-module
of rank $d$.  

Since $\ww\in V_L$, there
exist unique $\lambda_1,\ldots,\lambda_d\in L$ such that
\[ \ww = \lambda_1\vq_1+\cdots + \lambda_d\vq_d . \]
Also, since $\{\vq_1,\ldots,\vq_d\}$ is linearly independent over $L$, there exists
an $I\subseteq \{1,\ldots,n\}$ of order $d$ such that $0\neq \det(Q_I)\in S$.  For any
subset $I$ such that $\det(Q_I)\neq 0$, 
it follows that the vector $\left[ \lambda_1 \ \cdots \ \lambda_d\right]^T$ is the
unique solution to the equation $Q_I\vec{x} = (\ww)_I$.  By Cramer's rule we get
\[ \lambda_j = \frac{\mu_j}{\det(Q_I)}
 \ \ \mbox{ where } \ \mu_j = \det
\left[ \begin{array}{ccccc}  & | & | & | & \\
\cdots & (\vq_{j-1}))_I & (\ww)_I & (\vq_{j+1})_I & \cdots \\
& | & | & | & \end{array}\right] \]
for all $1\leq j\leq d$.  For each $j$, let $g_j=\gcd(\mu_j,\det(Q_I))$ and let us write
$\mu_j=h_jg_j$ and $\det(Q_I)=k_jg_j$ for some $h_j,k_j\in S$.  If $\mu_j\neq 0$, we
$\mu_j$ and $\det(Q_I)$ are homogenous and
$\deg(\mu_j)=\deg(\det(Q_I))=d$, hence $h_j$ and $k_j$ are
homogenous of the same degree of at most $d$.  If $\mu_j=0$, then $h_j=0$ and we let
$k_j=1$.  Finally, it is clear that $k_j$ must divide $\det(Q_I)$.  Since this is true for 
all $j$, we get $m={\mathrm{lcm}}(k_1,\ldots,k_d)$ must also divide $\det(Q_I)$ for 
all subsets $I\subseteq \{1,\ldots,n\}$ of order $d$, hence $\det(Q_I)\in \langle m\rangle$.  \qed

\begin{cor}\label{gcdcor}
Let $d<n$ and let $\{\vq_1,\ldots,\vq_d\}\subseteq (S_1)^n$ be linearly independent over the field
$L$.  Assume $\yy\in\spn_L(\vq_1,\ldots,\vq_d)$
and choose 
homogeneous elements
$h_1,\ldots,h_d,k_1,\ldots,k_d\in F[\yy]$ such that $k_1,\ldots,k_d$ are nonzero,
\[ \yy = \frac{h_1}{k_1}\,\vq_1 + \cdots + \frac{h_d}{k_d}\,\vq_d \]
and for each $1\leq j\leq d$ we have the pair $(h_j,k_j)$ are realtively prime of the same
total degree or if $h_j=0$ then $k_j=1$. Then $\deg(m) < d$ where 
$m$ is the least common multiple of $k_1,\ldots,k_d$.  Furthermore, if there exists an $I=\{i_1<\cdots<i_d\}|\subseteq \{1,\ldots,n\}$ such that 
$\det(Q_I)$ is irreducible, then $\yy\in\spn_F(\vq_1,\ldots,\vq_d)$.
\end{cor}

\prf   Assume the total degree of $m$ is $d$ and let $I\subseteq \{1,\ldots,n\}$ with $|I|=d$.  So either 
$\det(Q_I) = 0$ or
 the total degree of $\det(Q_I)$ is $d=\deg(m)$.  By Theorem \ref{yLthm}, $m|\det(Q_I)$ which implies $\det(Q_I) = \alpha m$ for some
 $\alpha \in F$.   
 Since $\yy\in \spn_{L}(\vq_1,\ldots,\vq_d)$, we know every $(d+1)\times (d+1)$ minor 
 of $C_{\yy}$ is zero.  Choose a subset $I'=\{i_1<\cdots<i_{d+1}\}\subseteq \{1,\ldots,n\}$ of order $d+1$ such that
 there is a subset $I\subset I'$ of order $d$ where $\det(Q_I)\neq 0$.  Therefore
 \[ y_{i_1}\det(Q_{I_1}) - y_{i_2}\det(Q_{I_2}) + \cdots + (-1)^dy_{i_{d+1}}\det(Q_{I_{d+1}}) = 0, \]
 where we let $I_t = \{ i_1,\ldots, i_{t-1},i_{t+1},\ldots, i_{d+1} \}$.  But
 by what we have done above, we know $\det(Q_{I_t}) = \alpha_t m$ for some 
 $\alpha_t\in F$ for all $1\leq t\leq d+1$.  Therefore, the above equation can be re-written as
 \[ m \left(\alpha_1y_{i_1} - \alpha_2y_{i_2} + \cdots + (-1)^d \alpha_{d+1}y_{i_{d+1}}\right)
 = 0. \]
 We know $m\neq 0$, so we must have
 \[ \alpha_1y_{i_1} - \alpha_2y_{i_2} + \cdots + (-1)^d\alpha_{d+1}y_{i_{d+1}} = 0. \]
 Since $\det(Q_I)\neq 0$, there exists at least one $\alpha_t\neq 0$, hence
\[ \alpha_1y_{i_1} - \alpha_2y_{i_2} + \cdots + (-1)^d\alpha_{d+1}y_{i_{d+1}} \neq 0, \]
which is a contradiction.  Hence, $\deg(m)<d$.

Finally, assume there exits an $I\subseteq \{1,\ldots,n\}$ such that $\det(Q_I)$ is irreducible.
Since $m|\det(Q_I)$ and $\deg(m)<d$, we get
$m\in F \Rightarrow k_j\in F$ for all $1\leq j\leq d$.  Therefore, for every $j$ either $h_j=0$ or 
$\deg(h_j)=\deg(k_j)=0$, so we get $\yy\in \spn_F(\vq_1,\ldots,\vq_d)$ as
claimed.  \qed

\vsp

Note that as an immediate application of Corollary \ref{gcdcor}, we see if $\yy\in\spn_L(\vq_1)$ for some $\vq_1\in (S_1)^n$, then $\yy\in\spn_F(\vq_1)$.  Hence for any $1$-dimensional $V\subseteq (S_1)^n$ with the $\yy$-local membership 
property must contain $\yy$.  Next we prove the analogous result when $V$ has
dimension $d=2$ as well.

\begin{thm}\label{d=1or2}
Let $n\geq 3$ and assume $F$ is algebraically closed and either $\chr(F)>n$ or $\chr(F)=0$.  Let $V\subseteq (S_1)^n$ be a two-dimensional $F$-subspace such that $S.V$ is a free $S$-module of 
rank $2$.  If $V$ has the $\yy$-local membership property, then 
$\yy\in V$.
\end{thm}

\prf  
Since $V$ has the $\yy$-local membership property, by Lemma \ref{Lspan}, 
we have $\yy\in V_L$. Let $\{\vq_1,\vq_2\}$ be an $F$-basis of $V$.  
So by Corollary \ref{gcdcor}, we can assume
\[ \yy = \frac{h_1}{k_1}\vq_1 + \frac{h_2}{k_2}\vq_2 \]
for some homogeneous polynomials $h_1,h_2,k_1,k_2\in S$ of degree at most 1 when $h_1\neq 0 \neq h_2$.  If $h_1=0$ or $h_2=0$, then by Corollary \ref{gcdcor}
 we get $\deg(k_1)=\deg(k_2)=0$, so
$\yy\in \spn_F(\vq_1,\vq_2)=V$. 
Therefore, we can assume $h_1\neq 0\neq h_2$ are of degree at most $1$
with $\deg(h_j)=\deg(k_j)$ for $j=1,2$.  

If $\deg(k_1)=\deg(k_2)=0$, we are done since $h_j/k_j\in F$ for $j=1,2$.  
So assume $\deg(k_1)=0$ and $\deg(k_2)=1$, so that $h_1/k_1=\alpha_1\in F$.  This implies 
\[ \frac{h_2}{k_2} \vq_2 = \yy - \alpha_1\vq_1 \Rightarrow k_2 | h_2q_{2,j} \mbox{ for all } 1\leq j\leq n. \]
Since $k_2$ is irreducible and $(h_2,k_2)$ are relatively prime, we get $k_2|q_{2,j}$ for all $1\leq j\leq n$, hence $\vq_2 = k_2\vv$ for some
nonzero $\vv\in F^n$.  Therefore, we have
\[ \yy = \alpha_1 \vq_1 + h_2\vv \Rightarrow \vq_1 = \frac{1}{\alpha_1} (\yy - h_2\vv) \]
for some $0\neq \alpha_1\in F$, $\vv\in F^n$, and $h_2\in S_1$.
Let ${\frak m}\subseteq S$ be a maximal ideal such that $k_2\in {\frak m}$ but $h_2\not\in {\frak m}$.   Then $\vq_2\in \ker(\phi_{{\frak m}})$.  
So, by our assumption, we have $\phi_{{\frak m}}(\yy)\in \spn_F(\phi_{{\frak m}}(\vq_1),
\phi_{{\frak m}}(\vq_2)) = \spn_F(\phi_{{\frak m}}(\vq_1))$.  Hence
there exists a nonzero $\beta\in F$ such
that $\alpha_1\phi_{{\frak m}}(\vq_1) + [h_2]_{\frak m}\vv = \phi_{{\frak m}}(\yy) = \beta \phi_{{\frak m}}(\vq_1)$.  So we get \\ 
$(\alpha_1-\beta)\phi_{{\frak m}}(\vq_1) + [h_2]_{\frak m}\vv = \vec{0}$.  
If $\alpha_1=\beta$, then $[h_2]_{\frak m}\vv = \vec{0}$ which implies $\vv=\vec{0}$ since $h_2\not\in {\frak m}$ 
which contradicts our choice of $\vq_2$.  If $\alpha_1\neq \beta$, then
we get ${\ds \phi_{{\frak m}}(\vq_1)-\frac{[h_2]_{\frak m}}{\beta-\alpha_1}\vv = \vec{0} \Rightarrow \phi_{{\frak m}}(\yy) - \left( \frac{\beta}{\beta-\alpha_1}\right)([h_2]_{\frak m})\vv 
= \vec{0}}$.  Hence $\phi_{{\frak m}}(\yy)\in \spn_F(\vv)$.  Now choose another maximal ideal ${\frak m}'\subseteq S$ containing $k_2$ but not $h_2$ with
such that $\phi_{{\frak m}}(\yy)$ and $\phi_{{\frak m}'}(\yy)$ are linearly independent in $F^n$.  (Note, we can do this as the 
Krull dimension of $S/\langle k_2\rangle$ is at least $2$.)  By the same argument, we get $\phi_{{\frak m}'}(\yy)\in \spn_F(\vv)$, 
which contradicts the choice of 
${\frak m}'$.

So we can assume $\deg(k_1)=\deg(k_2)=1$.  By Corollary \ref{gcdcor}, the degree of 
the lcm of $k_1$ and $k_2$ is at most 1, so 
we can assume $k_2 = \gamma k_1$ for some nonzero $\gamma\in F$.  By modifying our choice
of $h_2$ by a constant multiple if necessary, we can assume $k_1=k_2$
so that 
\[ \yy = \frac{1}{k_1} (h_1\vq_1 + h_2\vq_2)  \Rightarrow k_1\yy = h_1\vq_1 + h_2\vq_2.  \]
Let $\psi:S\rightarrow S/\langle h_2\rangle \cong F[w_1,\ldots,w_{n-1}]$ be the natural projection where we choose 
$\{ w_1,\ldots,w_{n-1}\}\subseteq \{\psi(y_1),\ldots,\psi(y_n)\}$ to be an $F$-linearly independent subset and let \[ \psi_n:S^n\rightarrow
S^n\otimes_{\psi}F[w_1,\ldots,w_{n-1}] \] be the induced map.  Hence we get 
\[ \psi(k_1)\psi_n(\yy) = \psi(h_1)\psi_n(\vq_1)  \Rightarrow \psi_n(\yy) = \frac{\psi(h_1)}{\psi(k_1)} \psi_n(\vq_1), \]
since $k_1\not\in\langle h_2\rangle$.  So using $S/\langle h_2\rangle$ instead of $S$ in Corollary 
\ref{gcdcor}, we get 
${\ds \frac{\psi(h_1)}{\psi(k_1)} \in F
\Rightarrow} $ $ \psi(h_1) = \gamma \psi(k_1)$ for some $\gamma\in F$, which gives us $h_1 = \gamma k_1 + \gamma' h_2$ for some
$\gamma'\in F$.  Therefore, we have
\[ \yy = \frac{h_1}{k_1}\vq_1 + \frac{h_2}{k_1}\vq_2 = \gamma \vq_1 + \frac{h_2}{k_1}(\gamma'\vq_1+\vq_2). \]
By replacing $\vq_2$ by $\gamma'\vq_1+\vq_2$ above, we reduce to the $\deg(h_1)=0$ case above, getting a contradiction.  \qed

\vsp

This establishes our local to global property when $d\leq 2<n$ and ${\mathrm{char}}(F)\neq 2,3$.   

We note that Theorem \ref{yLthm} and
Corollary \ref{gcdcor} only used our weaker hypothesis of assuming $\yy\in\spn_{L}(\vq_1,\ldots,\vq_d)$.  One can ask if there are stronger results if
we use the stronger hypothesis of $V=\spn(\vq_1,\ldots,\vq_d)$ having the $\yy$-local membership property.  
This is an area of future investigation, but the results in this section give a partial answer to when the 
$\yy$-local membership property holds in general.

For any $\{\vq_1,\ldots,\vq_{n-1}\}
\subseteq (S_1)^n$, we can re-write $C=C_{\yy}$ as
\[ C = y_1 A_1 + y_2 A_2 + \cdots + y_n A_n \]
for some matrices $A_1,\ldots,A_n\in M_n(F)$.  If $\{ \vq_1,\ldots,\vq_{n-1} \}$ is linearly
independent over $L$ and $\yy\in \spn_{L}(\vq_1,\ldots,\vq_{n-1})$, then we get 
$\det(C)=0$ in $S$.  Furthermore, by letting $y_j=1$ and
$y_i=0$ if $i\neq j$, it follows that 
$\det(A_j)=0$ and $A_j\vec{e}_n = \vec{e}_j$ for all $1\leq j\leq n$.

\begin{thm} \label{nontrivialint}
Let $\{\vq_1,\ldots,\vq_{n-1}\}\subseteq (S_1)^n$ be linearly independent over $L$.  Write
\[ C = C_{\yy}= y_1 A_1 + \cdots
+ y_nA_n \]
for some $A_1,\ldots,A_n\in M_n(F)$.  Then $\yy\in \spn_F(\vq_1,\ldots,\vq_{n-1})$ if
and only if $\{A_1,\ldots,A_n\}\subseteq I$ 
for some maximal left ideal $I\subseteq M_n(F)$.
\end{thm}  

\prf  First, let us assume that $\yy\in\spn_F(\vq_1,\ldots,\vq_{n-1})$.  Then there exists
$c_1,\ldots,c_{n-1}\in F$ such that 
\[ \vec{0} = c_1\vq_1+\cdots + c_{n-1}\vq_{n-1} + \yy = C \left[ \begin{array}{c}
c_1 \\ \vdots \\ c_{n-1} \\ 1 \end{array}\right]. \]
Let $I$ be the maximal left ideal ${\mathrm{l.ann}}\left(\left[ \begin{array}{c}
c_1 \\ \vdots \\ c_{n-1} \\ 1 \end{array}\right]\right)$ of $M_n(F)$.  
For arbitrary $j$, let $y_j=1$ and $y_i=0$ for $i\neq j$.  Then it follows that 
\[ A_j \left[ \begin{array}{c}
c_1 \\ \vdots \\ c_{n-1} \\ 1 \end{array}\right] = \vec{0} \ \Longrightarrow \  
\left[ \begin{array}{c}
c_1 \\ \vdots \\ c_{n-1} \\ 1 \end{array}\right] \in {\mathrm{Nul}}(A_j) \ \Longrightarrow 
\ A_j\in I. \]
Since $j$ was arbitrary, we get $\{A_1,\ldots,A_n\}\subseteq I$.

Conversely, assume $\{A_1,\ldots,A_n\}\subseteq I$ for some maximal 
left ideal $I$.  Then there
exists a nonzero $\vec{p}\in F^n$ such that $I={\mathrm{l.ann}}(\vec{p})$.  Hence $\vec{p} \in \bigcap_{j=1}^n {\mathrm{Nul}}(A_j)$.
Then 
\[ \vec{0} = C\vec{p} = p_1\vq_1 + \cdots + p_{n-1}\vq_{n-1} + p_n\yy \]
where $\vec{p} = \left[ \begin{array}{c} p_1 \\ \vdots \\ p_n\end{array}\right]$.  Note that
$p_n\neq 0$ since $p_n=0$ would imply that $\{\vq_1,\ldots,\vq_{n-1}\}$ is linearly
dependent over $F$, hence linearly dependent over $L$, contradicting our assumption.
Therefore, we get 
\[ \yy = \frac{-p_1}{p_n} \vq_1 + \cdots + \frac{-p_{n-1}}{p_n}\vq_{n-1} 
\in \spn_F(\vq_1,\ldots,\vq_{n-1}). \]
Hence our proposition follows.  \qed

\section{r1-free Subspaces}

As stated in our introduction, our interest in the above $\yy$-local membership property came about in our exploration of Mathieu-Zhao subspaces of matrix algebras
$M_n(F)$.  It follows from \cite[Thm.\ 4.2]{Zhao2} that $(0)\neq V\subseteq M_n(F)$
is a Mathieu-Zhao subspace of $M_n(F)$ if and only if $V$ does not contain any of the nonzero idempotents of $M_n(F)$.  In 2016, M. deBondt proved in \cite{deBondt} 
that any proper Mathieu subpace of $M_n(F)$ of codimension strictly
smaller than $n$ must be a subspace of the trace zero matrices of $M_n(F)$ when $F$ is algebraically closed and 
${\mathrm{char}}(F)=0$ or ${\mathrm{char}}(F)>n$.  In exploring deBondt's theorem further, we asked if we consider those
subspaces $V\subseteq M_n(F)$ that do not contain any rank $1$ idempotents, would an analogous result hold.  (We call such subspaces 
{\em r1-free subspaces}.)  In this section we will establish a correspondence between subspaces of $M_n(F)$ of codimension $d<n$ and
$d$-dimensional subspaces $V\subseteq (S_1)^n$ such that $S.V$ is a free $S$-submodule of rank $d$.  
Then $V$ defines a codimension $d$ 
r1-free subspace of $M_n(F)$ if and only if $V$ satisfies the $\yy$-local membership property.  
Furthermore, we get $\yy\in V$ if
and only if the corresponding codimension $d$ subspace of $M_n(F)$ is contained in the subspace of trace zero matrices.  
Hence Theorem \ref{d=1or2} tells us that deBondt's Theorem does not extend to r1-free subspaces unless you restrict to r1-free subspaces of codimension 1 or 2.

Recall that above we observed that for every $\vq\in (S_1)^n$, there exists a unique 
$b\in M_n(F)$ such that $\vq = b\,\yy$.  Therefore, for any subspace $V\subseteq (S_1)^n$
of dimension $d$,
we have a corresponding subspace 
\[ V^\flat = \{ b\in M_n(F) \ : \ \vq = b\,\yy \mbox{ for some } \vv\in V\}\subseteq M_n(F) \]
and the correspondence $V\mapsto V^\flat$ is a one-to-one correspondence of subspaces.
If we use the nondegenerate bilinear Killing form on $M_n(F)$ given by $\langle a,b\rangle =
{\mathrm{Tr}}(ab)$ for all $a,b\in M_n(F)$, we can consider the subspace 
$(V^\flat)^\perp \subseteq M_n(F)$
of codimension $d$ corresponding to $V$, hence we get a one-to-one correspondence between
$d$-dimensional subspaces $V$ of $(S_1)^n$  and codimension $d$ subspaces of $M_n(F)$.

In the masters thesis \cite{Kon} of A. Konijnenberg supervised by A. van den Essen, 
we have Theorem 3.7, which is a result attributed to deBondt.  Here we restate this theorem in terms of 
subspaces of $(S_1)^n$ having the
$\yy$-local membership property and include a proof for the reader's convenience.

\begin{thm}{\cite[Thm.\ 3.7]{Kon}} \label{Thm3.7}
Let $V\subseteq (S_1)^n$ have dimension $d<n$.  Then $V$ has the $\yy$-local membership
property if and only if $(V^\flat)^\perp$ is r1-free.
\end{thm}

\prf  First assume $V\subseteq (S_1)^d$ has the $\yy$-local membership property.  Let 
$e\in M_n(F)$ be a rank 1 idempotent in $(V^\flat)^\perp$.  
Then there exist $\uu,\vv\in F^n$ such that
$e=\uu\,\vv^{\,\tau}$ and $\vv^{\,\tau}\uu=1$, where we use the superscript $\tau$ to denote the transpose.  Hence for all
$b\in V^\flat$ we get $0={\mathrm{Tr}}(b\uu\,\vv^{\,\tau}) = \vv^{\,\tau}b\uu$.
Therefore, $\uu\not\in V^\flat\uu$.  But if we choose a maximal ideal ${\frak m}\subseteq S$ such
that $\phi_{\frak m}(\yy)=\uu$, we get $\uu\in \phi_{{\frak m}}(V)=V^\flat \uu$, which is a contradiction.  Therefore, $(V^\flat)^\perp$ is r1-free.

Conversely, let $V\subseteq (S_1)^n$ be such that $(V^\flat)^\perp$ is r1-free.  Let 
${\frak m}\subseteq S$ be a maximal ideal and let $\uu=\phi_{\frak m}(\yy)$.  If 
$\uu\not\in \phi_{\frak m}(V)=V^\flat\uu$, then there exists a $\vv\in F^n$ such that
$\vv^{\,\tau}\uu =1$ but $\vv^{\,\tau}\ww = 0$ for all $\ww\in \phi_{\frak m}(V)=V^\flat\uu$.  
Therefore, $e=\uu\,\vv^{\,\tau}$ is an idempotent in $(V^\flat)^\perp$, which contradicts our
assumption that $(V^\flat)^\perp$ is r1-free.  \qed

\vsp

So, combining Theorems \ref{d=1or2} and \ref{Thm3.7} with our Example in Section 2, we get
the following Theorem.

\begin{thm}\label{r1freecodim}
Assume $n\geq 3$ and
${\mathrm{char}}(F)=0$ or ${\mathrm{char}}(F)>n$.
Let $W\subseteq M_n(F)$ be an r1-free subspace of codimension $1$ or $2$.  Then $W$ is
a subspace of $I_n^\perp$, the subspace of trace zero matrices.  Furthermore, if $n\geq 4$ there exist
r1-free subspaces of codimension $3\leq d<n$ that are not subspaces of $I_n^\perp$.
\end{thm}

 
 \end{document}